\documentclass[12pt]{article}
\usepackage{amssymb, amsmath}
\usepackage{mathrsfs}       
\usepackage{color}  
     
\numberwithin{equation}{section}

\newtheorem{thm}{Theorem}[section]
\newtheorem{lem}{Lemma}[section]

\newtheorem{pro}{Proposition}[section]

\renewcommand{\Bbb}{\mathbb}

\newcommand{\e}{\varepsilon}
\renewcommand{\a}{a}

\renewcommand{\o}{\omega}

\newcommand{\bb}{\begin{equation}}
\newcommand{\ee}{\end{equation}}
\newcommand{\bq}{\begin{eqnarray}}
\newcommand{\eq}{\end{eqnarray}}
\newcommand{\bqn}{\begin{eqnarray*}}
\newcommand{\eqn}{\end{eqnarray*}}
\begin{document}
\title{Remarks on  type I blow up  for the 3D Euler equations and the 2D Boussinesq equations}
\author{Dongho Chae$^*$  and Peter Constantin$^\dagger$\\
\ \\
 $*$Department of Mathematics\\
Chung-Ang University\\
 Seoul 06974, Republic of Korea\\
 e-mail: dchae@cau.ac.kr\\
and \\
$\dagger$Department of Mathematics\\
Princeton University\\
Princeton, NJ 08544, USA\\
e-mail: const@math.princeton.edu }

\date{}
\maketitle

\begin{abstract}
In this paper we derive kinematic relations for quantities involving the rate of strain tensor and the Hessian of the pressure for solutions of the 3D Euler equations and the 2D Boussinesq equations.
Using these kinematic relations, we prove new blow up criteria and obtain conditions for the absence of type I singularity for these equations. We obtain both global and localized versions of the results.
Some of the new blow up criteria and type I conditions improve previous results of \cite{cc}.\\
\ \\
\noindent{\bf AMS Subject Classification Number:} 35Q31, 76B03\\
  \noindent{\bf
keywords:} Euler equations, Boussinesq equations, kinematic relations, blow up criterion, type I singularity
\end{abstract}

\section{The 3D Euler equations}
\setcounter{equation}{0}  
\subsection{Introduction and main results}

We consider the homogeneous incompressible  Euler equations on $\mathbb R^3\times [0, T)$.
$$
(E) \left\{\aligned  &u_t+u\cdot \nabla u= -\nabla p,\\
&\nabla \cdot u=0,\quad u(x,0)=u_0(x)
\endaligned \right.
$$
 where $u(x,t)=(u_1 (x,t), u_2 (x,t), u_3 (x,t))$ is the fluid velocity  and $p=p(x,t)$ is the  scalar pressure. 
We denote the initial velocity by $u_0 (x)= u(x,0) $  where $x\in \Bbb R^3$. For the Cauchy problem of the system (E) with $u_0 \in W^{2,q} (\Bbb R^3)$, $q>3$, $\nabla \cdot u_0 =0$, there exists a  local in time well-posedness  result \cite{kat}, but the question of the finite time blow up is a wide open problem. See e.g. \cite{maj, con1, con2, con3}  and the references therein for detailed discussions of the problem. For important partial results we refer \cite{bkm, cfm}. See also \cite{ker, luo} and references therein for  related numerical works. \\
 
 We associate to a solution $(u, p)$ of the Euler system (E) the  $\Bbb R^{3\times 3}$-valued functions 
$S=(S_{ij})$ and   $P=(P_{ij})$, where 
$$
S_{ij}= \frac12 (\partial_iu_j +\partial_j u_i ),\quad P_{ij}= \partial_i \partial_j p.
$$
   For the vorticity $\o=\nabla\times u$ we define the direction vectors
$$
 \xi= \o/|\o|, \qquad \zeta = S\xi/|S\xi |,
$$
 and  the scalar functions
$$
\alpha=S_{ij}\xi_i \xi_j, \quad   \rho= P_{ij}\xi_i\xi_j ,
$$
   where we used the convention of summing over repeated indices.  In the case $\o(x,t)=0$ we set $ \alpha(x,t)=\rho (x,t)= 0$. These quantities have been introduced previously  \cite{con3, maj,ch}.
 Note that $\xi $ is the {\em vorticity direction} vector, while $\zeta$ is  the {\em vorticity stretching direction} vector.  Below we also use the notations
  $
 [f]_+=\max\{ f, 0\}$ and $[f]_ - = \max\{- f, 0\}.$
 
  \begin{thm}
 Let $(u,p)\in C^1 (\mathbb R^3 \times (0, T))$ be a solution of the Euler equation (E) with $ u \in C([0, T); W^{2,q}  (\mathbb R^3))$, for some $ q> 3$.   Suppose the  following holds. Either 
 \begin{itemize}
\item[(i)] 
\bb\label{3a}
 \int_0 ^T\exp\left( \int_0^t \int_0 ^s\| [\zeta\cdot P\xi ]_- (\tau)  \|_{L^\infty} d\tau ds \right) dt <+\infty,
 \ee
 or
$$
 \int_0 ^T\exp\left( \int_0^t \int_0 ^s\| [|S\xi|^2 -2\alpha ^2 -\rho]_+(\tau)  \|_{L^\infty} d\tau ds\right) dt <+\infty,
$$
then $
 \limsup_{t\to T} \|u(t)\|_{ W^{2,q} } <+\infty.
$
 \item[(ii)]  If either 
\bb\label{3c}  \limsup_{t\to T} \,(T-t)^2 \| [\zeta\cdot P\xi ]_- (t)  \|_{L^\infty}  <1,
\ee
 or 
\bb\label{3d}
 \limsup_{t\to T} \,(T-t)^2 \| [|S\xi|^2 -2\alpha ^2 -\rho]_+(t)  \|_{L^\infty} <1,
\ee
 then $
\limsup_{t\to T} \|u(t)\|_{W^{2,q} } <+\infty.
$
 \end{itemize}
  \end{thm}
 {\bf Remark 1.1 } In \cite{cc} we obtained the above theorem with $ [\zeta\cdot P\xi   ]_- $ replaced by $ | P|$, which is the matrix norm of the Hessian of the pressure. Since $ | [\zeta\cdot P\xi   ]_- |\le |P|$ the above( and the localized version below) improve the results of Theorem 1.1 of \cite{cc}.  Furthermore, the above theorem implies  that
 the dynamical changes of the signs of the scalar quantities $\zeta\cdot P\xi$ and $|S\xi|^2 -2\alpha ^2 -\rho$ are important 
 in the phenomena of  blow up/regularity.\\
 \ \\

 The following is a localized version of the above theorem.
\begin{thm}
 Let $(u,p)\in C^1 (B(x_0, r) \times (T-r, T))$ be a  solution to (E)   with $ u \in C([T-r , T); W^{2, q}  ( B(x_0, r)))\cap L^\infty ( T-r, T; L^2 (B(x_0, r)))$ for some $q\in (3, \infty)$.  We suppose 
$$
 \int_{T-r} ^{T} \|u(t) \|_{L^\infty (B(x_0, r ))} dt<+\infty,
$$
and the following holds. Either
 \begin{itemize}
\item[(i)] 
$$
 \int_{T-r}^T\exp\left( \int_0^t \int_0 ^s\| [\zeta\cdot P\xi ]_-(\tau)  \|_{L^\infty (B(x_0, r))}  d\tau ds \right) dt <+\infty,
$$
 or
$$
 \int_{T-r}^T\exp\left( \int_{T-r}^t \int_{T-r} ^s\| [|S\xi|^2 -2\alpha ^2 -\rho]_+(\tau)   \|_{L^\infty (B(x_0, r))} d\tau ds\right) dt <+\infty.
$$
Then for all $\e\in (0,r)$ $
\limsup_{t\nearrow T} \|u(t)\|_{ W^{2,q} (B(x_0,\e))} <+\infty.
$

 \item[(ii)]  If  \eqref{3c} holds, and  if either 
\bb\label{3g} \limsup_{t\to T} \,(T-t)^2 \| [\zeta\cdot P\xi ]_-(t)  \|_{L^\infty (B(x_0, r))}  <1,
\ee
 or 
\bb\label{3h}
 \limsup_{t\to T} \,(T-t)^2 \| [|S\xi|^2 -2\alpha ^2 -\rho]_+(t)   \|_{L^\infty (B(x_0, r))} <1,
\ee
 then for all $\e\in (0,r)$ $
\limsup_{t\nearrow T} \|u(t)\|_{ W^{2,q} (B(x_0,\e))} <+\infty.
$
 \end{itemize}
\end{thm}

\subsection{Kinematic relations}

We use  the particle trajectory mapping  $\a \mapsto  X ({\a} , t)$ from $\mathbb R^3$ into $\mathbb R^3$   generated by $u=u(x,t)$, which means the solution of  the ordinary differential equation,
$$
\left\{ \aligned &\frac{\partial   X (\a, t)}{\partial t}   =u( X (\a, t)  , t) \quad \text{on}\quad (0,T),\\
                    & X (\a , 0)  =\a \in \mathbb R^3.\endaligned
                    \right.
$$
The material derivative of $f=f(x,t)$ is defined by 
$$D_t f:=\partial_t f +u\cdot \nabla f.
$$
We note that
$
(D_t f) (X(a, t),t)=\frac{\partial }{ \partial t } \left\{ f (X(\a, t),t)  \right\}.
$
\begin{pro} 
Let $(u,p)$ be a solution of (E), which belongs to $C^1(\Bbb R^3 \times (0, T) )$. We use the above notations.
Then, the followings hold true on $\Bbb R^3 \times (0, T)$.
\begin{align} \label{5}
D_t |S\o |&= - \zeta \cdot P\o, \\
\label{6}
D_t ^2\log |\o| &= |S\xi|^2 -2\alpha^2  -\rho , \\
\label{7}
  (D_t  |\o|)^2+   (|D_t \xi ||\o|)^2 & = |S\o|^2, \\
 \label{8}
  (D_t  |S\o|)^2+   (|D_t \zeta ||S\o|)^2 &= |P\o|^2,\\
  \label{8a}
         \end{align} 
\end{pro}
 {\bf Remark 1.2 } Applying the inequality, $ a_1 +\cdots +a_n\le \sqrt{ n(a_1 ^2 +\cdots +a_n^2)} $ to equations \eqref{7}, \eqref{8} and \eqref{9} respectively, we obtain the following differential inequalities with the coefficients consisting of derivatives of the direction fields $\xi$ and  $\zeta$,
 \begin{align}
 D_t |\o| + |D_t \xi| |\o|&\le \sqrt2|S\o|,\label{8b}\\
 D_t |S\o| + |D_t \zeta ||S\o|&\le \sqrt2 |P\o|, \label{8c}\\
 D_t |S\o| +  |D_t \zeta| D_t |\o| +  |D_t \zeta| |D_t \xi | |\o| &\le \sqrt3 |P\o|.\label{8d}
 \end{align}
 For an  implication  of \eqref{8b}  combined with \eqref{8c}, in particular,  see  Remark 1.3 below.\\
 \ \\
  {\bf Proof  of Proposition 1.1 }
Taking the gradient of (E), we find 
\bb\label{9}
D_t \nabla u= -(\nabla u)^2 -P.
\ee
We observe the decomposition of the matrix,
 $$\nabla u= S +\Omega, \quad \text{where} \quad \Omega_{ij}= \frac12 (\partial_i u_j-\partial_j u_i) =  \frac12 \epsilon_{ijk} \o_k.
 $$
Here,  $\epsilon_{ijk}$ is the  totally skew-symmetric tensor with normalization $\epsilon_{123}=1$.
Taking  the skew symmetric part of \eqref{9}, we obtain  the vorticity equations
   \bb\label{10}
    D_t  \o = \o \cdot \nabla u =S\o,
\ee
where we used the fact 
$
\o_j \partial_j u_i =\o_j S_{ji} +\frac12\epsilon_{jik}  \o_j \o_k =\o_j S_{ji}.
$
Taking the symmetric part of \eqref{9}, on the other hand,  we find
\bb\label{11} 
D_t S=-S^2 +\frac14 ( |\o|^2 I -\o\otimes \o )-P.
\ee
Contracting  \eqref{10}  with $\o$, and dividing the both sides by $|\o|^2$, we have 
\bb\label{12} D_t  |\o|  = \alpha |\o|.
\ee
 From \eqref{10} and \eqref{12}   we derive 
\bb\label{13}
D_t \xi = \frac{D_t \o}{|\o|} - \o \frac{D_t  |\o|}{|\o|^2} =      S \xi-\alpha \xi.
\ee
Applying $D_t$ to  \eqref{10},  using \eqref{11}, we find
\begin{align}\label{14}
D_t ^2 \o&= (D_t S) \o +S D_t \o= -S^2 \o -P \o + S^2 \o \cr
&= -P\o,
\end{align} 
which was the key kinematic relation used  in \cite{cc}.
Multiplying \eqref{14}  by $D_t \o=S\o$ from the left, we obtain
$$
|D_t \o|  D_t  |D_t \o|= \frac12 D_t |D_t \o|^2 =D_t \o \cdot D_t^2 \o= -S\o \cdot P\o,
$$
Dividing the both sides by $ |D_t \o|= |S\o|$, we find
\bb\label{14a} 
D_t |D_t \o |=D_t |S\o|= -\zeta\cdot P\o,
\ee
and \eqref{5} is proved.
Now we prove \eqref{6}. Observing  $\xi \cdot D_t  \xi=0$,   we compute
 \begin{align}\label{15}
 D_t ^2 |\o|&= D_t \{ \xi \cdot D_t  (|\o| \xi )\}= D_t \xi \cdot  D_t \o  + \xi \cdot  D_t ^2 \o\cr
 & =(S -\alpha I ) \xi \cdot  S\o   -\xi \cdot  P\o  =\left(|S\xi|^2- \alpha^2  - \rho\right) |\o|.
 \end{align}
We divide \eqref{15} by $|\o|$, then using \eqref{12}, we deduce
 \begin{align*}
 & |S\xi|^2- \alpha^2  - \rho=\frac{D_t ^2  |\o|}{ |\o|} = D_t \left( \frac{D_t  |\o|} {|\o| }\right) + \frac { (D_t  |\o|)^2 }{|\o|^2} 
 = D_t ^2 \log |\o|+  \alpha^2. 
 \end{align*} 
 The formula \eqref{6} is proved.
 Taking the square of  \eqref{13}, and multiplying it by $|\o|^2$,  we have
    \bb\label{17}
 |S\o|^2= \alpha^2|\o|^2 +|D_t \xi|^2 |\o|^2  =( D_t |\o| )^2 + |D_t \xi|^2 |\o|^2,
 \ee  
 and, \eqref{7} is proved.
To show \eqref{8} we compute, using \eqref{14} and \eqref{14a}, 
\begin{align*}
 D_t \zeta &= \frac{ D_t ^2  \omega}{ | D_t  \omega|} - \frac{ D_t \omega  D_t \left(  |D_t \omega|\right)}{ | D_t  \omega|^2}=-\frac{ P \omega}{ |S\omega|}+\frac{ S  \omega  (\zeta \cdot P \xi ) |\omega|}{  |S  \omega |^2}\cr
 &=\frac{ -P\xi  +  (\zeta \cdot P \xi ) \zeta}{ |S\xi|}.
  \end{align*}
Because $D_t\zeta $ is perependicular to $\zeta$, in view of the fact that $\zeta $ has unit length, this yields  an orthogonal decomposition of $P\xi$,
 \bb\label{18a}
 P\xi=  (\zeta \cdot P \xi ) \zeta- |S\xi | D_t \zeta= (\zeta \cdot P \xi ) \zeta +  \frac{(D_t \zeta \cdot P \xi )}{|D_t \zeta |^2}  D_t \zeta,  \ee 
 which implies
 \bb\label{19}
\frac{D_t \zeta}{|D_t \zeta|}\cdot P \xi= -|S\xi| |D_t \zeta|.
 \ee
 The decomposition \eqref{18a}, combined with \eqref{14a} and \eqref{19},  implies by the Pythagoras theorem
  \begin{align}\label{20} 
 |P\o|^2 &=    (\zeta \cdot P\o )^2+   \left(\frac{D_t \zeta}{|D_t \zeta|} \cdot P\o\right)^2  \cr
 &= (D_t  |S\o|)^2+   |D_t \zeta |^2|S\o|^2.
  \end{align}
The inequality \eqref{8} follows from this immediately. 
Substituting \eqref{17} into \eqref{20}, we have \eqref{8a}. $\square$\\
\ \\

\subsection{ Proofs of the main theorems}
In order to prove Theorem 1.1 we shall use the following lemma.
\begin{lem}
 Let $\alpha =\alpha (t)$ be a non-decreasing  function, and $\beta =\beta(t)\ge 0$ on $[a,b]$.
 \begin{itemize}
 \item[(i)] Suppose $y=y(t)$ satisfies
 $$
 y(t) \le \alpha (t) + \int_a ^t \beta (\tau ) y(\tau ) d\tau \quad \forall t\in [a, b].
 $$
 Then, for all $t\in (a, b]$ we have
$$
y(t)\le  \alpha (t)  \exp \left( \int_a ^t \beta (\tau)d\tau \right) .
$$
 \item[(ii)] We assume furthermore $ y(t) \ge 0$ on $[a, b]$.  Suppose 
 $$
 y(t) \le \alpha (t) + \int_a ^t \int_a ^s \beta (\tau ) y(\tau )d\tau ds \quad \forall t\in [a, b].
 $$
 Then, for all $t\in (a, b]$ we have
$$
y(t)\le  \alpha (t)  \exp \left( \int_a ^t \int_a ^s \beta (\tau)d\tau ds \right) .
$$
\end{itemize}
 \end{lem}

 {\bf Proof }   In the case (i)  from the well-known Gronwall inequality and the assumption of  non-decreasing property of $\alpha$ we have
\begin{align*} 
 y(t) &\le \alpha (t)+ \int_a ^t \alpha (s) \beta (s) \exp \left(\int_s ^t  \beta (\tau ) d\tau \right) ds \le \alpha (t)+ \alpha (t) \int_a ^t \beta (s) \exp \left(\int_s ^t  \beta (\tau ) d\tau \right) ds  \cr
 &=\alpha (t)-  \alpha (t) \int_a ^t  \frac{d}{ds} \left\{ \exp \left(\int_s ^t  \beta (\tau ) d\tau \right)\right\} ds= \alpha (t)  \exp \left( \int_a ^t \beta (\tau)d\tau \right).
 \end{align*} 
For the case (ii) we  observe
$$
 y(t)\le \alpha (t) + \int_a ^t \int_0 ^s \beta (\tau) y(\tau) d\tau ds
 \le \alpha (t) + \int_a ^t  \sup_{ a< \tau< s}  y(\tau)\int_a ^s  \beta (\tau ) d\tau ds.
$$
 Since the function $t\mapsto  \alpha (t) + \int_a ^t  \sup_{ a< \tau< s} y(\tau) \int_a ^s  \beta (\tau ) d\tau ds$ is  non-decreasing  on $[a, b]$,  setting $ h(t) =\int_a^ t \beta (s)ds$ and $Y(t)= \sup_{ a< \tau< t}y(\tau)$, we have
 $$ Y(t) \le \alpha (t)+ \int_a ^t Y(s) h(s) ds.$$
 Applying (i), we obatin
$$
 y(t)\le Y(t)\le \alpha (t) \exp\left( \int_a^t h (s )ds\right)=  \alpha (t) \exp\left( \int_a^t \int_a ^s \beta(\tau) d\tau ds\right)
$$
 for all $t\in [a, b]$. $\square$\\
 \ \\

\noindent{\bf Proof of  Theorem 1.1 and Theorem 1.2  } 
We integrate  \eqref{5} along the trajectory  for $t\in [0, s]$ to find
  \begin{align*} 
&\frac{\partial}{\partial s} | \o (X(a, s) ,s) | \le   \left|\frac{\partial}{\partial s}\o (X(a, s),s)\right| = | (D_s \o ) (X(a,s),s)| = |S\o (X(a, s),s) |\cr
&\quad= |S_0 (a)\o _0 (a)| -\int_0 ^s (\zeta \cdot P \xi )(X(\a, \tau), \tau ) |\o (X (a, \tau), \tau |d\tau,
\end{align*}
 from which, after integrating with respect to $s$ over $[0, t]$,  we have 
$$
 |\o (X(\a, t),t) |\le |\o_0 (\a) | + |S_0 (\a)\o _0 (\a)|t +
 \int_0 ^t \int_0 ^s[\zeta \cdot P \xi]_- (X(\a, \tau), \tau ) |\o (X (a, \tau), \tau)|  d\tau  ds.
$$
 Applying  Lemma 2.1(ii) to solve this differential inequality,  we find 
\begin{align}\label{26}
 |\o (X(\a, t),t) |&\le ( |\o_0 (\a) | + |S_0 (\a)\o _0 (\a)|t ) \times \cr
 &\qquad \times\exp\left(\int_0 ^t \int_0 ^s  [\zeta \cdot P\xi]_- (X(\a, \tau), \tau ) d\tau ds \right).
\end{align}
Taking the supremum over $a\in \Bbb R^3$, and integrating it with respect to $t$ over  $[0, T]$, we find
\begin{align}\label{27}
\int_0 ^T \|\o(t)\|_{L^\infty} dt &\le ( \|\o_0 \|_{L^\infty} + \|S_0 \o_0 \|_{L^\infty} T)\times \cr
 &\qquad \times\int_0 ^T \exp\left(\int_0 ^t \int_0 ^s  \|[\zeta \cdot P \xi ]_-(\tau )\|_{L^\infty} d\tau ds \right)dt.
\end{align}
Integrating \eqref{6}  twice with respect to the time variable over $ [0, s]$,  we hvae
\bb\label{28}
|\o (X(a,t), t) |= |\o_0 (a) | \exp \left( \int_0 ^t \int_0 ^s [ |S\xi|^2 -2\alpha^2  -\rho]_+ (X(a, \tau), \tau ) d\tau ds \right),
\ee
and therefore
$$
\int_0 ^T \|\o(t)\|_{L^\infty} dt\le  \|\o_0 \|_{L^\infty} \int_0 ^T  \exp \left( \int_0 ^t \int_0 ^s\| [|S\xi|^2 -2\alpha^2  -\rho]_+ (\tau) \|_{L^\infty} d\tau ds \right)dt.
$$
Applying the well-known Beale-Kato-Majda criterion \cite{bkm}  to \eqref{27} and \eqref{28}, we obtain the desired conclusion of 
Theorem 1.1(i). 
The argument of proof of  Theorem 1.1(ii), using the result of (i) is the similar to \cite{cc}, and we ommit it here. \\
\ \\   
The proof of Theorem 1.2, using  the key pointwise estimates of the vorticity along the trajectories, \eqref{26}  and \eqref{27} is similar to the corresponding 
ones in \cite{cc}, and we do not repeat it here.
$\square$\\
\ \\
{\bf Remark 1.3  }   The linear differential inequalities \eqref{8b} and \eqref{8c}  along the trajectory can be solved as
    \begin{align}\label{30}
   & |\o(X(a, t),t)| \le |\o_0(a)| e^{- \int_0 ^t  |D_t \xi  ( X(a,s),s) | ds } \cr
    &\qquad+\sqrt2 \int_0 ^t |S\o (X(a,s), s)|e^{-\int_s ^t | D_\tau\xi (X(a,\tau), \tau)| d\tau } ds,
    \end{align}
    and 
      \begin{align}\label{31}
   & |S\o(X(a, t),t)| \le |S_0\o_0(a)| e^{- \int_0 ^t  |D_t \zeta ( X(a,s),s) | ds }   \cr
   &\qquad  +\sqrt2 \int_0 ^t |P\xi  (X(a,s), s)|  |\o (X(a,s), s)| e^{-\int_s ^t | D_\tau\zeta (X(a,\tau), \tau)| d\tau } ds    \end{align}
    respectively. Parenthetically one can also use \eqref{8a} to deduce  
    $$D_t  |S\o|   +|D_t \zeta| |D_t \xi | |\o| \le  \sqrt2  |P\o |, $$
    and then
        \begin{align*}
   & |S\o(X(a, t),t)|  \le  |S_0\o_0(a)| e^{- \int_0 ^t  |D_t \zeta ( X(a,s),s) |  |D_t \xi( X(a,s),s) | ds }   \cr
   & +\sqrt2 \int_0 ^t |P\xi  (X(a,s), s)|  |\o (X(a,s), s)| e^{-\int_s ^t | D_\tau\zeta (X(a,\tau), \tau)|  | D_\tau\xi (X(a,\tau), \tau)|d\tau } ds
   \end{align*}
    instead of \eqref{31}.     Inserting \eqref{31} into   \eqref{30}, we find
   \begin{align}\label{32}
  & |\o(X(a, t),t)| \le |\o_0 (a)| + \sqrt2 |S_0 (a) \o_0 (a) | \int_0 ^t  e^{- \int_0 ^s |D_\tau  \zeta| d\tau  }    e^{- \int_s ^t |D_\tau  \xi | d\tau  }  ds\cr  
    &+2 \int_0 ^t \int_0 ^s  |P\xi  ( X(a, \sigma), \sigma ) | |\o ( X(a, \sigma), \sigma )|e^{-\int_{\sigma} ^s | D_\tau  \zeta| d\tau}  e^{-\int_s ^t  |D_\tau\xi | d \tau} d\sigma   ds
       \end{align} 
   Applying Lemma 1.1 to \eqref{32}, we obtain
   \begin{align}\label{33}
   & |\o(X(a, t),t)|  \le   \left( |\o_0 (a)| + \sqrt2 |S_0 (a) \o_0 (a) | t \right) \times \cr
   &\times \exp\left(  2 \int_0 ^t \int_0 ^s  |P\xi  ( X(a, \sigma), \sigma ) | e^{-\int_{\sigma} ^s |  D_\tau\zeta (X(a,\tau), \tau)|  d\tau}  e^{-\int_s ^t  | D_\tau\xi (X(a,\tau), \tau) | d \tau} d\sigma   ds  \right). \cr
   \end{align} 
   Since the quantities expressing the magnitudes of the changes of the two direction  vectors    $D_t \xi$ and $D_t \zeta$   contribute to the integral in the right hand side of \eqref{33}  through factors like  $e^{-\int_s ^t  | D_\tau\xi (X(a,\tau), \tau) | d \tau}$, they appear to have a desingularizing effect for the vorticity. We do not know, however a way to exploit this effect in the blow up criterion and the absence of the type I blow up.
   If we ignore the factor $e^{-\int_s ^t  | D_\tau\xi (X(a,\tau), \tau) | d \tau}$   in \eqref{33}, 
   taking supremum over $a\in \Bbb R^3$, and integrating over $t\in [0, T]$  then we have  an estimate
\begin{align*} 
\int_0 ^T \|\o(t)\|_{L^\infty} dt&\le \left( \|\o_0 \|_{L^\infty} + \sqrt2 \|S_0 \o_0  \|_{L^\infty}  T\right)\times \cr
&\qquad \times \int_0 ^T \exp\left(2\int_0 ^t\int_0 ^s
 \|P\xi (\tau)  \|_{L^\infty} d\tau ds\right) dt
 \end{align*}
which yields a blow up criterion  weaker than \eqref{3a}.

\section{The 2D Boussinesq equations}
\setcounter{equation}{0}  
Here we are concerned with the homogeneous incompressible  Boussinesq equation on $\Bbb R^2$.
$$
(B) \left\{\aligned  &u_t+u\cdot \nabla u= -\nabla p + \theta e_2,\label{e1}\\
& \theta_t + u\cdot \nabla \theta=0,\\
&\nabla \cdot u=0,
\endaligned \right.
$$
 where $u(x,t)=(u_1 (x,t), u_2 (x,t))$ is the fluid velocity  and $p=p(x,t)$ is the   pressure, and
 $\theta=\theta(x,t)$ is the temperature. Let $u_0(x)= u(x,0), \theta_0 (x,0)$ be the initial data of the system (B).
The local well-posedness for the Boussinesq system for $(u_0, \theta_0 )\in W^{2,q} (\Bbb R^2)$, $q>2$,  is well-known(see e.g.\cite{cn}), but the 
 question of finite time blow up is a wide open problem similarly to the case of the 3D Euler equations. It is also well-known that there exists a strong similarity between (B) and the axisymmetric solution of the 3D Euler equations(see e.g.\cite{maj}). \\
 \ \\
 For a solution $(u, p, \theta)$ of the system (B) let us introduce the  $\Bbb R^{2\times 2}$-valued functions 
$U=( \partial_iu_j)$ and   $P=(\partial_i \partial_j p)$. For the vector filed $\nabla ^\perp \theta=(-\partial_2\theta, \partial_1 \theta)$
 we define the direction vectors
$$
 \xi= \nabla ^\perp \theta/|\nabla ^\perp \theta|, \qquad \zeta = U\nabla ^\perp \theta/|U\nabla ^\perp \theta |,
$$
 and  the scalar functions
$$
\alpha=\xi\cdot U\xi, \quad   \rho= \xi\cdot P\xi .
$$
 
  \begin{thm}
 Let $(u,p)\in C^1 (\mathbb R^2 \times (0, T))$ be a solution of the Boussinesq  equation (B) with $ u \in C([0, T); W^{2,q}  (\mathbb R^2))$, for some $ q> 2$.   Suppose the  following holds. Either 
 \begin{itemize}
\item[(i)] 
\bb\label{35}
 \int_0 ^T (T-t) \exp\left( \int_0^t \int_0 ^s\| [\zeta\cdot P\xi ]_- (\tau) \|_{L^\infty} d\tau ds \right) dt <+\infty,
 \ee
 or
 \bb\label{36}
 \int_0 ^T(T-t) \exp\left( \int_0^t \int_0 ^s\| [|U\xi|^2 -2\alpha ^2 -\rho]_+(\tau)   \|_{L^\infty} d\tau ds\right) dt <+\infty,
 \ee
then $
 \limsup_{t\to T} \|u(t)\|_{ W^{2,q} } <+\infty.
$
 \item[(ii)]  If either 
\bb\label{36a}
 \limsup_{t\to T} \,(T-t)^2 \| [\zeta\cdot P\xi ]_- (t) \|_{L^\infty}  <2,
\ee
 or 
\bb\label{36b}
 \limsup_{t\to T} \,(T-t)^2 \| [|U\xi|^2 -2\alpha ^2 -\rho]_+(t)   \|_{L^\infty} <2,
\ee
 then $
\limsup_{t\to T} \|u(t)\|_{W^{2,q} } <+\infty.
$
 \end{itemize}
  \end{thm}
 {\bf Remark 2.1 } Note the relaxed smallness condition of for the  nonexistence of type I blow up in \eqref{36a} and  \eqref{36b} compared to \eqref{3c} and \eqref{3d}  respectively  in the case of 3D Euler equations. This is due to the extra factor,
  $(T-t)$ in   \eqref{35} and \eqref{36}, which originate from 
  the non blow up criterion, $\int_0 ^T (T-t) \|\nabla ^\perp \theta  (t)\|_{L^\infty}  dt<+\infty$  in  \cite[Theorem 1.2 (ii)]{cw}. \\  
  \ \\
 The following is a localized version of the above theorem.
\begin{thm}
 Let $(u,p)\in C^1 (B(x_0, r) \times (T-r, T))$ be a  solution to (E)   with $ u \in C([T-r , T); W^{2, q}  ( B(x_0, r)))\cap L^\infty ( T-r, T; L^2 (B(x_0, r)))$ for some $q\in (2, \infty)$. Let us assume
$$
 \int_{T-r} ^{T} \|u(t) \|_{L^\infty (B(x_0, r ))} dt<+\infty.
$$
If either
 \begin{itemize}
\item[(i)] 
\bb\label{38}
 \int_{T-r}^T(T-t) \exp\left( \int_0^t \int_0 ^s\| [\zeta\cdot P\xi ]_- (\tau)  \|_{L^\infty (B(x_0, r))}  d\tau ds \right) dt <+\infty,
 \ee
 or
 \bb\label{39}
 \int_{T-r}^T (T-t) \exp\left( \int_{T-r}^t \int_{T-r} ^s\| [|U\xi|^2 -2\alpha ^2 -\rho]_+(\tau)  \|_{L^\infty (B(x_0, r))} d\tau ds\right) dt <+\infty,
 \ee
then for all $\e\in (0,r)$ $
\limsup_{t\nearrow T} \|u(t)\|_{ W^{2,q} (B(x_0,\e))} <+\infty.
$

 \item[(ii)]  If either 
\bb
\label{39a}
  \limsup_{t\to T} \,(T-t)^2 \| [\zeta\cdot P\xi ]_- (t)  \|_{L^\infty (B(x_0, r))}  <2,
\ee
 or 
\bb
\label{39b}
 \limsup_{t\to T} \,(T-t)^2 \| [|U\xi|^2 -2\alpha ^2 -\rho]_+(t)   \|_{L^\infty (B(x_0, r))} <2,
\ee
 then for all $\e\in (0,r)$ $
\limsup_{t\nearrow T} \|u(t)\|_{ W^{2,q} (B(x_0,\e))} <+\infty.
$
 \end{itemize}
\end{thm} 
 {\bf Remark 2.2 } Similarly to Remark 2.1 we  also note here   relaxed smallness condition of for the  nonexistence of type I blow up in \eqref{39a} and  \eqref{39b} compared to \eqref{3g} and \eqref{3h}  respectively. This is due to the extra factor,
  $(T-t)$ in  \eqref{38} and \eqref{39}, which are from  the local version of the  non blow up criterion,  $\int_0 ^T (T-t) \|\nabla^\perp \theta (t)\|_{L^\infty (B(x_0, r))}  dt<+\infty $ in \cite[Theorem 2.1]{cw1}.

\subsection{Kinematic relations}
\begin{pro} 
Let $(u,p, \theta)$ be a solution of (B), which belongs to $C^1(\Bbb R^2 \times (0, T) )$. We use the above notations.
Then, the followings hold true on $\Bbb R^2 \times (0, T)$.
\begin{align} \label{40}
D_t |U\nabla^\perp \theta |&= - \zeta \cdot P\nabla^\perp \theta, \\
\label{41}
D_t ^2\log |\nabla^\perp \theta| &= |U\xi|^2 -2\alpha^2  -\rho , \\
\label{42}
  (D_t  |\nabla^\perp \theta|)^2+   (|D_t \xi ||\nabla^\perp \theta|)^2 & = |U\nabla^\perp \theta|^2, \\
 \label{43}
  (D_t  |U\nabla^\perp \theta|)^2+   (|D_t \zeta ||U\nabla^\perp \theta|)^2 &= |P\nabla^\perp \theta|^2,\\
  \label{44}
  (D_t  |U\nabla^\perp \theta|)^2+ ( |D_t \zeta| D_t |\nabla^\perp \theta|)^2   +(|D_t \zeta| |D_t \xi | |\nabla^\perp \theta| )^2&=  |P\nabla^\perp \theta |^2.  \end{align} 
\end{pro}
{\bf Remark  2.1 } Although the above results look similar to those in Proposition 1.1 we have essentially different features because we do not use the symmetric part of $U$ and because there exists no relation between $\nabla^\perp \theta $ and the skew symmetric part of $U$. \\
\ \\

{\bf Proof of Proposition 2.1  } Taking $\nabla$ on the first equation of $(B)$, we find
 \bb\label{45}
 D_t U + U^2 =-P  + \nabla (\theta e_2 ),
 \ee
 Taking $\nabla^\perp $ on the second  equation of $(B)$, 
 \bb\label{46}
 D_t \nabla^\perp \theta= U \nabla^\perp \theta  .
 \ee  
Let us compute
 \begin{align}\label{47}
 D_t ^2   \nabla^\perp \theta &=   D_t U \nabla^\perp \theta  + U D_t  \nabla^\perp \theta  \cr
 &=  -  U^2 \nabla^\perp \theta  - P\nabla^\perp \theta   +  U^2 \nabla^\perp \theta   +
 \nabla^\perp \theta \cdot  \nabla (\theta e_2 ) \cr
  &= - P\nabla^\perp \theta,
\end{align}
where we used the fact 
\bb\label{47a}
\nabla^\perp \theta \cdot  \nabla (\theta e_2 ) =0.
\ee
We multiply \eqref{47} by $ D_t \nabla^\perp\theta$ to have
\begin{align*}
 |D_t \nabla^\perp \theta| D_t  |D_t \nabla^\perp \theta|&=  \frac12D_t \left(  |D_t\nabla^\perp \theta |^2\right)  = D_t \nabla^\perp \theta \cdot D_t ^2 \nabla^\perp \theta\cr
 & = -U\nabla^\perp \theta \cdot P \nabla^\perp \theta.
 \end{align*} 
Dividing the both sides by $ |D_t \nabla^\perp \theta|= |U \nabla^\perp \theta|$, we find
$$
 D_t |D_t \nabla^\perp \theta|= D_t |U\nabla^\perp \theta| =  -\zeta \cdot P \xi |\nabla^\perp \theta|,
$$
 and \eqref{40} is proved.
Multiplying  \eqref{46}  by $\nabla ^\perp \theta$, we deduce
 \bb\label{48}
  D_t |  \nabla^\perp \theta|=\alpha |  \nabla^\perp \theta|.
  \ee
  Using \eqref{46} and \eqref{48}, we compute
\bb\label{49}
 D_t \xi= \frac{ D_t \nabla^\perp \theta}{ | \nabla^\perp \theta|} - \frac{ \nabla^\perp \theta D_t | \nabla^\perp \theta|}{ | \nabla^\perp \theta|^2}=  U\xi -\alpha \xi.
\ee
 This can be viewed  as an orthogonal decomposition of $U\xi$,
$$
 U\xi= \alpha \xi + D_t\xi= \alpha \xi + \frac{ D_t \xi \cdot U\xi }{ |D_t \xi |^2 } D_t \xi ,
$$
 which  shows
 \bb\label{50}
 D_t \xi \cdot U\xi = |D_t \xi|^2= |U\xi|^2-\alpha^2= |U\xi|^2-\frac{ (D_t | \nabla^\perp \theta|)^2}{| \nabla^\perp \theta|^2}.
 \ee 
 Multiplying the both sides of \eqref{50} by $| \nabla^\perp \theta|^2$, the formula \eqref{42} follows immediately.
 Using  \eqref{45} and \eqref{49},  we compute
 \begin{align*}
  &D_t ^2 \log |  \nabla^\perp \theta|  =D_t  \alpha =D_t (\xi\cdot U\xi )\cr
  &\qquad= D_t \xi\cdot  U \xi + \xi \cdot D_t U \xi + \xi \cdot U D_t \xi\cr
  &\qquad= (U\xi -\alpha \xi) \cdot U\xi + \xi \cdot (-U^2 -P)\xi +  \xi \cdot U (U\xi-\alpha \xi )\cr 
  &\qquad= |U\xi|^2-2 \alpha ^2 -\rho,
     \end{align*} 
    where we used $ \xi\cdot \nabla (\theta e_2)=0$, which follows from \eqref{47a}. The formula \eqref{41} is proved.
Using \eqref{47} and \eqref{40}, we compute
\begin{align}\label{52} 
 D_t \zeta &= \frac{ D_t ^2  \nabla^\perp \theta}{ | D_t  \nabla^\perp \theta|} - \frac{ D_t \nabla^\perp \theta  D_t \left(  |D_t \nabla^\perp \theta|\right)}{ | D_t  \nabla^\perp \theta|^2}\cr
 &=-\frac{ P \nabla^\perp \theta}{ |U\nabla^\perp \theta|}+\frac{ U  \nabla^\perp \theta  (\zeta \cdot P \xi ) |\nabla^\perp \theta|}{  |U  \nabla^\perp \theta |^2}\cr
 &=\frac{ -P\xi  +  (\zeta \cdot P \xi ) \zeta}{ |U\xi|}.
  \end{align}
 The formula \eqref{52}  yields  an orthogonal decomposition of $P\xi$,
 \bb\label{53}
 P\xi=  (\zeta \cdot P \xi ) \zeta- |U\xi | D_t \zeta= (\zeta \cdot P \xi ) \zeta +  \frac{(D_t \zeta \cdot P \xi )}{|D_t \zeta |^2}  D_t \zeta,  \ee 
 which implies
 \bb\label{54}
\frac{D_t \zeta}{|D_t \zeta|}\cdot P \xi= -|U\xi| |D_t \zeta|.
 \ee
The decomposition \eqref{53} also implies by the Pythagoras theorem, and then using \eqref{40} and \eqref{54}, 
 \begin{align*}
 |P\nabla^\perp \theta|^2 &=    (\zeta \cdot P\nabla^\perp \theta  )^2+   \left(\frac{D_t \zeta}{|D_t \zeta|} \cdot P \nabla^\perp \theta \right)^2  \cr
 &=   (D_t  |U\nabla^\perp \theta|)^2+ |D_t \zeta|^2 |U\nabla^\perp \theta |^2,
 \end{align*}
      which verifies \eqref{43}. Inserting the expression of $|U\nabla^\perp \theta|^2$  in \eqref{43} into \eqref{44}, we obtain \eqref{45}. 
      $\square$
 
\subsection{Proof of the main results}
\noindent{\bf Proof of  Theorem 2.1 and Theorem 2.2  }  The proof is  similar to the case of 3D Euler equations. The main difference is that here we start from the kinematic relations of  the Boussinesq equations in Proposition 2.1.
Integrating  \eqref{40} along the trajectory  for $t\in [0, s]$, we obtain
  \begin{align*} 
&\frac{\partial}{\partial s} |\nabla^\perp \theta (X(a, s) ,s) | \le   \left|\frac{\partial}{\partial s}\nabla^\perp \theta(X(a, s),s)\right| = | (D_s \nabla^\perp \theta ) (X(a,s),s)| = |U\nabla^\perp \theta (X(a, s),s) |\cr
&\quad= |S_0 (a)\o _0 (a)| -\int_0 ^s( \zeta \cdot P \xi )(X(\a, \tau), \tau ) |\o (X (a, \tau), \tau |d\tau.
\end{align*}
After integrating this again with respect to $s$ over $[0, t]$,  we find 
\begin{align*}
 |\nabla^\perp \theta (X(\a, t),t) |&\le |\nabla^\perp \theta_0 (\a) | + |\nabla^\perp \theta_0 (a) \cdot \nabla u_0 (\a)|t \cr
 &\qquad+
 \int_0 ^t \int_0 ^s[\zeta \cdot P \xi ]_-(X(\a, \tau), \tau ) |\nabla^\perp \theta (X (a, \tau), \tau)|  d\tau   ds.
\end{align*} 
 Thanks to Lemma 2.1(ii)   we find 
\begin{align}\label{56}
 |\nabla^\perp \theta(X(\a, t),t) |&\le ( |\nabla^\perp \theta_0 (\a) | + |\nabla^\perp \theta_0 \cdot \nabla u_0  (\a)(\a)|t )\times \cr
&\qquad \times  \exp\left(\int_0 ^t \int_0 ^s  [\zeta \cdot P \xi ]_-(X(\a, \tau), \tau ) d\tau ds \right).
\end{align}
Taking the supremum over $a\in \Bbb R^2$, and integrating it with respect to $t$ over  $[0, T]$ after multiplying by $T-t$, we find
\begin{align}\label{57}
&\int_0 ^T(T-t)  \|\nabla^\perp \theta(t)\|_{L^\infty} dt \le ( \|\nabla^\perp \theta_0 (\a \|_{L^\infty} + \|\nabla^\perp \theta_0 \cdot \nabla u_0 \|_{L^\infty} T) \times\cr
&\qquad\times\int_0 ^T (T-t)\exp\left(\int_0 ^t \int_0 ^s  \|[\zeta \cdot P \xi (\tau )]_- \|_{L^\infty} d\tau ds \right)dt.
\end{align}
Integrating \eqref{41}  twice with respect to the time variable over $ [0, s]$,  we hvae
\bb\label{58}
|\nabla^\perp \theta (X(a,t), t) |\le |\nabla^\perp \theta_0 (a) | \exp \left( \int_0 ^t \int_0 ^s  [|U\xi|^2 -2\alpha^2  -\rho]_+ (X(a, \tau), \tau ) d\tau ds \right),
\ee
and from which we also deduce 
\begin{align}\label{59}
&\int_0 ^T (T-t)\|\nabla^\perp \theta(t)\|_{L^\infty} dt\le  \|\nabla^\perp \theta_0 \|_{L^\infty} \times \cr
&\qquad \times\int_0 ^T (T-t) \exp \left( \int_0 ^t \int_0 ^s\| [ |U\xi|^2 -2\alpha^2  -\rho]_+ (\tau) \|_{L^\infty} d\tau ds \right)dt.
\end{align}
Applying the blow up criterion  of \cite[Theorem 1.2 (ii)]{cw} to \eqref{57} and \eqref{59}, we obtain the desired conclusion of 
Theorem 2.1(i). 
The proof of  Theorem 2.1(ii), using the result of (i) is the similar to the one in \cite{cc}.\\
\ \\   
The proof of Theorem 2.2, using  the key estimates \eqref{56}  and \eqref{58} is also similar to the corresponding ones in \cite{cc}. The essential point here is that we apply the local version of the blow up criterion \cite[Theorem 2.1]{cw1}. $\square$\\
\ \\

      $$ \text{\bf Acknowledgements} $$
Chae was partially supported by NRF grant 2021R1A2C1003234, while the work of Constantin was partially supported by the Simons Center for Hidden Symmetries and Fusion Energy. 

  $$ \text{\bf  Conflicts of Interest Statement} $$
The authors declare  that there is no conflict of interest.

\end{document}